\documentclass[11pt,a4paper]{amsart}
\usepackage[utf8]{inputenc}
\usepackage[T1]{fontenc}
\usepackage[english]{babel}

\usepackage{sidecap}

\usepackage{verbatim}
\usepackage{lmodern}
\usepackage{amsmath}
\usepackage{amssymb}
\usepackage{amsthm}
\usepackage{thmtools}
\usepackage{enumitem}
\usepackage{graphicx}
\usepackage[hidelinks]{hyperref}
\usepackage{mathtools}
\usepackage{indentfirst}
\usepackage{tabularx}
\usepackage{bbm}
\usepackage[capitalise]{cleveref}
\usepackage{stmaryrd}
\usepackage{tabularx}

\usepackage[left=1.2in,top=2cm,right=1.2in,bottom=2cm]{geometry}
\usepackage{array}
\newcolumntype{C}[1]{>{\centering\let\newline\\\arraybackslash\hspace{0pt}}m{#1}}

\usepackage[
textwidth=2cm,
textsize=small,
colorinlistoftodos]
{todonotes} 

\usepackage{tikz} 
\usepackage{tikz-cd} 
\newsavebox{\InnerDiagram}
\newsavebox{\RightDiagram}
\usepackage{quiver}
\usepackage{adjustbox}
\newcommand{\btk}{\begin{tikzcd}}
\newcommand{\etk}{\end{tikzcd}}
\usetikzlibrary{arrows,calc,decorations.markings}
\usepackage{xparse}

\usepackage{blindtext}

\makeatletter
\newcommand{\@bbify}[1]{
  \ifcsname b#1\endcsname
  \message{WARNING: Overwriting b#1 with blackboard letter!}
  \fi
  \expandafter\edef\csname b#1\endcsname
  {\noexpand\ensuremath{\noexpand\mathbb #1}\noexpand\xspace}}
\newcommand{\@calify}[1]{
  \ifcsname c#1\endcsname
  \message{WARNING: Overwriting c#1 with calligraphic letter!}
  \fi
  \expandafter\edef\csname c#1\endcsname
  {\noexpand\ensuremath{\noexpand\mathcal #1}\noexpand\xspace}}
\newcommand{\@bfify}[1]{
  \ifcsname bf#1\endcsname
  \message{WARNING: Overwriting c#1 with bold letter!}
  \fi
  \expandafter\edef\csname bf#1\endcsname
  {\noexpand\ensuremath{\noexpand\mathbf #1}\noexpand\xspace}}
\newcounter{@letter}\stepcounter{@letter}
\loop\@bbify{\Alph{@letter}}\@calify{\Alph{@letter}}\@bfify{\Alph{@letter}}
\ifnum\the@letter<26\stepcounter{@letter}\repeat
\makeatother

\tikzstyle{d}=[double distance=.3ex]
\tikzstyle{w}=[preaction={draw=white,-,line width=5pt}]

\newcounter{diagram}
\renewcommand{\thediagram}{\thetheorem}

\tikzset{%
node distance=1.5cm, la/.style={scale=0.8}, lasmall/.style={scale=0.75}, over/.style={auto=false,fill=white,inner sep=1.5pt, minimum size=0, outer sep=0},
    symbol/.style={%
        draw=none,
        every to/.append style={%
            edge node={node [sloped, allow upside down, auto=false]{$#1$}}},

    }, pro/.style={postaction={decorate,decoration={
        markings,
        mark=at position .5 with {\node at (0,0) {$\bullet$};}
      }},
      inner sep=.9ex,
      },
      prosmall/.style={postaction={decorate,decoration={
        markings,
        mark=at position .5 with {\node at (0,0) {$\scriptstyle \bullet$};}
      }},
      inner sep=.9ex,
      },
  n/.style={double equal sign distance, -implies}, t/.style={double distance=2.5pt, -implies, postaction={draw,-}},
}

\tikzset{%
node distance=1.5cm, la/.style={scale=0.8}, rr/.style={xshift=1.5cm},
space/.style={xshift=.5cm},
    symbol/.style={%
        draw=none,
        every to/.append style={%
            edge node={node [sloped, allow upside down, auto=false]{$#1$}}},

    }
}

\newcommand{\Set}{\mathrm{Set}}

\newcommand{\id}{\mathrm{id}}
\newcommand{\ob}{\mathrm{ob}}
\newcommand{\op}{\mathrm{op}}

\newcommand{\cat}{\mathrm{Cat}}

\DeclareMathOperator{\hocolim}{\mathrm{hocolim}}
\DeclareMathOperator{\diag}{\mathrm{diag}}

\newlist{rome}{enumerate}{7}
\setlist[rome]{label=(\roman*)}

\newtheorem{theorem*}{Theorem}
\newtheorem{theorem}{Theorem}[section]

\declaretheorem[name=Theorem,numbered=yes]{theoremA}

\theoremstyle{definition}
\newtheorem{defn}[theorem]{Definition}

\newtheorem{ex}[theorem]{Example}
\newtheorem{notation}[theorem]{Notation}

\theoremstyle{remark}
\newtheorem{rem}[theorem]{Remark}

\crefname{theorem}{Theorem}{Theorems}
\crefname{cor}{Corollary}{Corollaries}
\crefname{prop}{Proposition}{Propositions}
\crefname{lem}{Lemma}{Lemmas}

\crefname{defn}{Definition}{Definitions}
\crefname{terminology}{Terminology}{Terminologies}
\crefname{ex}{Example}{Examples}
\crefname{notation}{Notation}{Notations}
\crefname{descr}{Description}{Descriptions}
\crefname{constr}{Construction}{Constructions}

\crefname{rem}{Remark}{Remarks}

\usepackage{blindtext}

\renewcommand\thepart{\Roman{part}.}

\makeatletter
\renewcommand\part{%
  \par
  \addvspace{4ex}%
  \@afterindenttrue
  \secdef\@part\@spart
}

\def\@part[#1]#2{%
    \ifnum \c@secnumdepth >\m@ne
      \refstepcounter{part}%

      \addcontentsline{toc}{section}{\hspace{-.5cm} \bfseries\thepart\hspace{1em}#1}%
    \else
      \addcontentsline{toc}{section}{#1}%
    \fi
    {\parindent \z@ \raggedright
     \interlinepenalty \@M
     \normalfont
     \thispagestyle{empty}
     \ifnum \c@secnumdepth >\m@ne
      \centering\large\textsc{\textbf{\thepart}}\nobreakspace
     \fi
     \centering\large\textsc{\textbf{#2}}
     \par}%
    \nobreak
    \vskip .3cm
    \@afterheading}
\def\@spart#1{%
      \addcontentsline{toc}{part}{#1}%
    {\parindent \z@ \raggedright
     \interlinepenalty \@M
     \normalfont
     \thispagestyle{plain}
     \centering\large\textsc{\textbf{#1}}
     \par}%
    \nobreak
    \vskip .3cm
    \@afterheading}
\makeatother 

\title{Thomason's colimit theorem for the double category of elements}

\author[A.\ Gill]{Andrew Gill}
\address{School of Mathematics, University of Minnesota, Minneapolis MN, 55415, USA}
\email{gill0565@umn.edu}

\author[M.\ Sarazola]{Maru Sarazola}
\address{School of Mathematics, University of Minnesota, Minneapolis MN, 55415, USA}
\email{maru@umn.edu}

\begin{document}

\begin{abstract}
We show that, for any 2-category $\cC$ and 2-functor $F\colon\cC^\op\to\underline{\cat}$, the double category of elements $\iint_\cC F$ introduced by Grandis and Par\'e satisfies a version of Thomason's colimit theorem; that is, there is a weak homotopy equivalence $B\hocolim F\simeq B(\iint_\cC F)$.
\end{abstract}

\maketitle

\section{Introduction}

The category of elements construction
\[ \textstyle \int_\cC \colon \cat(\cC^{\op},\Set)\to \cat_{/\cC}\]
gives an equivalence between the categories of functors from a fixed category $\cC^\op$ to $\Set$, and of discrete fibrations over $\cC$.
It is intimately linked with the study of representable functors, as a classical result shows that a functor $F\colon \cC^{\op}\to \Set$ is representable if and only if its category of elements $\int_\cC F$ has a terminal object; see e.g.~\cite[Proposition 2.4.8]{Riehlcontext}. Hence, the category of elements gives us a way to characterize representable functors, and through them, universal properties, which are then used to understand key constructions such as adjunctions and (co)limits.

While this explains the central place that the category of elements occupies in category theory, it turns out that this construction---together with its non-discrete analogue, the Grothendieck construction---is also quite relevant in homotopy theory, mainly due to the following theorem by Thomason.
\begin{theoremA}\cite[Theorem 1.2]{thomason}
   For any category $\cC^\op$ and functor $F\colon\cC^\op\to \cat$, there is a weak homotopy equivalence \[B\hocolim  F\simeq B\left(\int_\cC F\right).\]
\end{theoremA}
Put in words, this result says that the Grothendieck construction of a diagram of categories gives a categorical model for the homotopy type of the homotopy colimit of the diagram of categories. Of course, the particular instance of this theorem for the discrete case of a functor $F\colon\cC^\op \to \Set$ also holds.

Interestingly, one can extend the category of elements to the 2-categorical setting. Now categories are enriched over $\cat$ instead of $\Set$, and so a 2-categorical version must take as input a 2-functor $F\colon\cC^\op\to\underline{\cat}$ from a 2-category $\cC^\op$ to the 2-category $\underline{\cat}$, and produce as output a 2-category $\int_\cC F$. Such a construction has been discussed by several authors, and can be found for instance as the ``discrete'' version of \cite[\S 4]{colims2functors} or \cite[Construction 2.2.1]{Buckley} as well as the underlying strict version of \cite[Definition 4.1]{cegarra}, \cite[\S 4]{grothcat1}, \cite[\S 3]{groth2}, among others.

From a homotopy-theoretic point of view, this 2-categorical construction enjoys the same useful properties as its 1-dimensional analogue; that is, it admits a version of Thomason's colimit theorem.

\begin{theoremA}\cite[Theorem 4.5(i)]{cegarra}\label{thm:2catofelements}
    For any 2-category $\cC^\op$ and 2-functor $F\colon\cC^\op\to \underline{\cat}$, there is a weak homotopy equivalence \[B\hocolim  F\simeq B\left(\int_\cC F\right).\]
\end{theoremA}

From a purely categorical point of view, its behavior is somewhat mixed. On the one hand, this construction gives the expected biequivalence between the 2-category of 2-functors $2\text{-}\cat(\cC^\op,\underline{\cat})$ and a 2-category of discrete 2-fibrations over $\cC$, as shown in \cite[Theorem 3.15]{lambert}. On the other hand, this 2-category of elements fails to capture the desired universal properties; that is, one cannot encode the representability of a 2-functor in terms of 2-terminal objects in its 2-category of elements. A modern account of this fact can be found in \cite{cM1}---where the authors thoroughly show that, in fact, various pseudo and lax versions fail as well.

The solution, suggested by Verity \cite{VerityThesis} and Grandis--Par\' e \cite{GrandisPare}, is to pass to the setting of double categories. In \cite[\S 1.2]{GraParPersistentII}, Grandis--Par\'e introduce a construction that takes as input a 2-functor $F\colon\cC^\op\to\underline{\cat}$ for some 2-category $\cC^\op$, and returns a \emph{double category of elements} $\iint_\cC F$. This construction now has all the desired categorical features: it encodes the representability of 2-functors $F\colon\cC^\op\to\underline{\cat}$ in terms of the existence of double terminal objects in the double category of elements $\iint_\cC F$ \cite[Theorem 6.8]{cM2}; allows us to express weighted 2-limits in terms of (conical) double limits \cite[Proposition~1.4]{GraParPersistentII}; and gives an equivalence between the category of 2-functors $2\text{-}\cat(\cC^\op,\underline{\cat})$ and a category of internal discrete fibrations over $\mathbb{H}\cC$ \cite[Theorem 5.1]{MSVenrichedGC}.

Given the suitability of this double category of elements to replace its 2-categorical predecessor for categorical purposes, the goal of this article is to explore whether this construction still produces the correct topological space. We give a positive answer in our main result, where we prove that the double category of elements also satisfies a version of Thomason's colimit theorem.

\begin{theoremA}[\cref{thm:main}]
    For any 2-category $\cC^\op$ and 2-functor $F\colon\cC^\op\to \underline{\cat}$, there is a weak homotopy equivalence \[B\hocolim  F\simeq B\left(\iint_\cC F\right).\]
\end{theoremA}

\subsection*{Conventions} Throughout the paper, all notions of 2-categories, 2-functors, and double categories are strict. We use $\underline{\cat}$ to denote the 2-category of categories, functors, and natural transformations. Given a 2-category $\cC$, we use $\cC^\op$ to denote the 2-category in which the direction of the morphisms is reversed, but not that of the 2-cells. We use $\alpha\vert\beta$ to denote the horizontal pasting of 2-cells or of squares as below
\[\begin{tikzcd}
    \bullet\rar[bend right=40]\rar[bend left=40]\rar[phantom, "\Downarrow\scriptstyle{\alpha}" description] & \bullet\rar[bend right=40]\rar[bend left=40]\rar[phantom, "\Downarrow\scriptstyle{\beta}" description] & \bullet
\end{tikzcd} \qquad \qquad
\begin{tikzcd}
    \bullet\dar\rar\ar[dr,phantom,"\Downarrow\scriptstyle{\alpha}" description]& \bullet\dar\rar\ar[dr,phantom,"\Downarrow\scriptstyle{\beta}" description] & \bullet\dar\\
  \bullet\rar & \bullet\rar  & \bullet
\end{tikzcd}\] and $\beta\circ\alpha$ to denote the vertical composite of 2-cells of the form below.
\[\begin{tikzcd}
    \bullet\rar\rar[phantom, bend right=35,"\Downarrow\scriptstyle{\beta}" description]\rar[phantom, bend left=35,"\Downarrow\scriptstyle{\alpha}" description]\rar[bend right=75]\rar[bend left=75] & \bullet
\end{tikzcd}\]
Finally, for any object $x$, horizontal morphism $f$, and vertical morphism $\varphi$ in a double category, we use $\id_x$ to denote the corresponding horizontal identity, $e_x$ for the vertical identity, $e_f$ for the vertical identity square on $f$, and $\id_\varphi$ for the horizontal identity square on $\varphi$.

\section{Background}

In this section we provide the reader with some of the necessary background by describing the main characters in our theorem: the double category of elements and the homotopy colimit of a 2-functor $F\colon\cC^\op\to\underline{\cat}$. We assume a previous familiarity with the basics of 2-categories, double categories, and (bi)simplicial objects, and refer the readers looking for introductory accounts on these notions to \cite{johnsonyau,grandis,goerssjardine}, respectively.

\subsection{Categories of elements for 2-functors}

We start by introducing the relevant categories of elements for a 2-functor $\cC^\op\to\underline{\cat}$, first defining the 2-categorical version and then the double categorical one.

\begin{defn}\label{defn:2catelt}
Let $\cC$ be a 2-category and $F\colon\cC^\op\to\underline{\cat}$ be a 2-functor. The \textbf{2-category of elements} $\int_\cC F$ consists of the following data:
\begin{itemize}
    \item its objects are pairs $(c,x)$ where $c$ is an object of $\cC$ and $x$ is an object of $Fc$;
    \item its morphisms are pairs $(f,\varphi)\colon (c,x)\to (c',x')$ given by a morphism $f\colon c\to c'$ in $\cC$ and a morphism $\varphi\colon x\to Ff(x')$ in $Fc$;
    \item its 2-cells $(f,\varphi)\Rightarrow (g,\psi)$ are 2-cells $\alpha\colon f\Rightarrow g$ in $\cC$ such that the following diagram commutes.
    \[\begin{tikzcd}
        x\rar["\varphi"]\ar[dr,"\psi"'] & Ff(x')\dar["(F\alpha)_{x'}"]\\
        & Fg(x')
    \end{tikzcd}\]
\end{itemize}
The identity morphism at $(c,x)$ is $(\id_c, \id_x)$, and the identity 2-cell at $(f,\varphi)$ is $\id_f$. The composite of morphisms \[(c,x)\xrightarrow{(f,\varphi)} (c',x')\xrightarrow{(g,\psi)} (c'',x'')\] is given by the pair $(gf, Ff(\psi)\circ\varphi)$, and composition of 2-cells works as in $\cC$.
\end{defn}

\begin{defn}\label{defn:dblcatelt}
   Let $\cC$ be a 2-category and $F\colon\cC^\op\to\underline{\cat}$ be a 2-functor. The \textbf{double category of elements} $\iint_\cC F$ consists of the following data:
        \begin{itemize}
            \item its objects are pairs $(c,x)$ where $c$ is an object of $\cC$ and $x$ is an object of $Fc$;
            \item its horizontal morphisms $(c,x)\to(c',x')$ are morphisms $f\colon c\to c'$ in $\cC$ such that $Ff(x')=x$;
            \item its vertical morphisms $(c,x)\to (c',x')$ only exist when $c=c'$ and consist of morphisms $\varphi\colon x\to x'$ in $Fc$;
            \item its squares
            \[\begin{tikzcd}
            (c,Ff(x'))\rar["f"]\dar["\psi"']\ar[rd,phantom, "\Downarrow\scriptstyle{\alpha}" description] & (c',x')\dar["\varphi"]\\
            (c, Fg(y'))\rar["g"'] & (c', y')
            \end{tikzcd}\] are 2-cells $\alpha\colon f\Rightarrow g$ in $\cC$ such that $\psi=(F\alpha)_{y'}\circ Ff(\varphi)$.
        \end{itemize}
        Identities and composites of horizontal morphisms work as in $\cC$, and those of vertical morphisms work as in $Fc$. Horizontal and vertical composition of squares, as well as identities, work as in $\cC$.
\end{defn}

We now exhibit a small example that illustrates both the differences and the similarities in the two structures defined above.

\begin{ex}\label{ex1}
    Let $\cC$ be the 2-category given by the following data
    \begin{center}
        \begin{tikzpicture}
            \node at(-0.2,0) {$c$};
            \node at(1.2,0.05) {$c'$};
            \node at(0.52,0) {$\Downarrow\scriptstyle{\alpha}$};
            \draw[->] (0, 0.2) to[bend left] node[above] {$\scriptstyle{f}$} (1,0.2);
       \draw[->] (0, -0.2) to[bend right] node[below] {$\scriptstyle{g}$} (1,-0.2);
        \end{tikzpicture}
    \end{center}
    and consider the 2-functor $F\colon\cC^\op\to\underline{\cat}$ defined as below. 

\begin{lrbox}{\InnerDiagram}
        \begin{tikzcd}[row sep=tiny]
            & Ff(x')\ar[dd, "\alpha"]\\
            x\ar[ur]\ar[dr]\\
            & Fg(x')
        \end{tikzcd}
    \end{lrbox}

    \begin{center}
    \begin{tikzpicture}
        \node at(-1.75,0) {$\Downarrow\scriptstyle{\alpha}$};
         \node at(0,0) {$\{x'\}=Fc'$};
         \node at(-5,0) {$\begin{tikzcd}
        Fc= \left\{ \usebox{\InnerDiagram} \right\}
    \end{tikzcd}$};
       \draw[->] (-1, 0.2) to[bend right] node[above] {$\scriptstyle{Ff}$} (-2.5,0.2);
       \draw[->] (-1, -0.2) to[bend left] node[below] {$\scriptstyle{Fg}$} (-2.5,-0.2);
    \end{tikzpicture}
    \end{center}

The data in the 2-category of elements $\int_\cC F$ can be represented in the following two diagrams:
    \begin{equation}\label{eq1}
    \begin{tikzcd}[column sep=small]
        & (c,Ff(x'))\ar[rd,"{(f,\id)}"]\\
        (c,x)\ar[rr, phantom, bend left=30, "\Downarrow\scriptstyle{\id}"]\ar[rr, phantom, bend right=30, "\Downarrow\scriptstyle{\id}"]\ar[ur,"{(\id,\varphi)}"]\ar[dr,"{(\id,\psi)}"']\ar[rr, bend left=15, "{(f,\varphi)}" description]\ar[rr, bend right=15, "{(g,\psi)}" description]\ar[rr, phantom, "\Downarrow" pos=0.45, "\scriptstyle{\alpha}"' pos=0.55] & & (c',x')\\
        & (c, Fg(x'))\ar[ur,"{(g,\id)}"']
    \end{tikzcd} \qquad \qquad
    \begin{tikzcd}[column sep=small]
         & (c,Ff(x'))\ar[rd,"{(f,\id)}"]\ar[dd,"{(\id,(F\alpha)_{x'})}" description]\\
        (c,x)\ar[rr,phantom, "\Downarrow\scriptstyle{\id}" pos=0.12]\ar[rr,phantom, "\Downarrow\scriptstyle{\alpha}" pos=0.87]\ar[ur,"{(\id,\varphi)}"]\ar[dr,"{(\id,\psi)}"'] & & (c',x')\\
        & (c, Fg(x'))\ar[ur,"{(g,\id)}"']
    \end{tikzcd}
    \end{equation}
 On the other hand, the data in the double category of elements $\iint_\cC F$ is given as follows:
    \begin{equation}\label{eq2}
    \begin{tikzcd}
       (c,x)\rar[equal]\ar[dd,"\psi"'] \ar[ddr,phantom, "\Downarrow\scriptstyle{\id}" description]& (c,x)\dar["\varphi"']\\
       &  (c,Ff(x'))\rar["f"]\dar["(F\alpha)_{x'}"']\ar[dr,phantom, "\Downarrow\scriptstyle{\alpha}" description] & (c',x')\dar[equal]\\
      (c,Fg(x'))\rar[equal] &  (c,Fg(x'))\rar["g"'] & (c,x')
    \end{tikzcd}\end{equation}

    The similarities between both sets of data are apparent: they have the same objects, and in both instances we are recording in some form the information of $f,g,\varphi,\psi$, and $\alpha$. But we can also note some important differences:
    \begin{itemize}
        \item In the 2-category of elements, we have a morphism $(f,\varphi)$ which is the composite of the morphisms $(\id,\varphi)$ and $(f,\id)$. However, in the double category of elements, the analogue of $(\id,\varphi)$ is given by the vertical morphism $\varphi$, and the analogue of $(f,\id)$ is the horizontal morphism $f$; thus, there is no notion of composite of these different types of maps. 
        \item In the double category of elements, there is a constraint in the objects that may participate in a non-identity square: the second coordinate in the objects on the left must be of the form $Ff(x')$ for some appropriate $f$ and $x'$. In contrast, the 2-category of elements has non-trivial 2-cells involving objects of general form, such as the 2-cell $\alpha$ on the left in \cref{eq1}.
        \item The 2-cell $\alpha\in\cC$ gives rise to two 2-cells in the 2-category of elements, of the form $(f,\varphi)\Rightarrow (g,\psi)$ and $(f,\id)\Rightarrow (g, (F\alpha)_{x'})$. On the other hand, it gives rise to a unique square in the double category of elements as pictured in \cref{eq2}.
    \end{itemize}
\end{ex}

\subsection{Homotopy colimit of a 2-functor}

The homotopy colimit of a 2-functor was defined in \cite[Definition 2.2.2]{hocolim}, with the purpose of studying the algebraic $K$-theory groups of an exact category.  Curiously, its definition will not be instrumental in the proof of our main theorem, but for the sake of completeness and to acknowledge the major role it plays in our result, we recall its definition.

\begin{defn}
    Let $\cC$ be a 2-category and $F\colon\cC^\op\to\underline{\cat}$ be a 2-functor. The \textbf{homotopy colimit} of $F$ is the simplicial category \[\hocolim F\colon\Delta^\op\to\cat\] whose category of $n$-simplices is given by
    \[(\hocolim F)_n=\bigsqcup_{(c_0,\dots,c_n)\in\ob\cC^{n+1}} Fc_0\times \cC^\op (c_0,c_1)\times\dots\times \cC^\op (c_{n-1},c_n).\]
    The face functor $d_0$ maps a component of the form $Fc_0\times \cC^\op (c_0,c_1)\times\dots\times \cC^\op (c_{n-1},c_n)$ to $Fc_1\times \cC^\op (c_1,c_2)\times\dots\times \cC^\op (c_{n-1},c_n)$ and is induced by the functor \[Fc_0\times \cC^\op (c_0,c_1)\to F c_1\] which acts on objects and morphisms as
    \begin{center}
    \begin{tikzpicture}
    \node at(1.4,1.2) {$(x, c_0\xleftarrow{f} c_1) \mapsto Ff(x)$};
            \node at(-1,0) {$(x\xrightarrow{\varphi} y, c_0$};
            \node at(3.8,0) {$c_1)\mapsto Ff(x)\xrightarrow{Fg(\varphi)\circ(F\alpha)_x} Fg(y)$};
            \node at(0.52,-0.1) {$\Downarrow\scriptstyle{\alpha}$};
            \draw[->] (1, 0.1) to[bend right] node[above] {$\scriptstyle{f}$} (0,0.1);
       \draw[->] (1, -0.3) to[bend left] node[below] {$\scriptstyle{g}$} (0,-0.3);
        \end{tikzpicture}
        \end{center}
    All other face maps, and all degeneracy maps, act only on the hom-categories as the corresponding nerve maps.
\end{defn}

\begin{rem}
The homotopy colimit of a 2-functor yields a topological space, for instance\footnote{There are several ways to produce a space out of a simplicial category, all of which are homotopy equivalent, and hence indistinguishable for our purposes.} by considering the bisimplicial set whose set of $(m,n)$-simplices is given by $N_m (\hocolim F)_n$, and taking the geometric realization of its diagonal simplicial set $[k]\mapsto N_k (\hocolim F)_k$. We denote this space by $B\hocolim F$.
\end{rem}

\section{Thomason's colimit theorem}

The goal of this section is to prove our main theorem, which we restate here.

\begin{theorem}\label{thm:main}
 For any 2-category $\cC^\op$ and 2-functor $F\colon\cC^\op\to \underline{\cat}$, there is a weak homotopy equivalence \[B\hocolim  F\simeq B\left(\iint_\cC F\right).\]
\end{theorem}

In light of \cref{thm:2catofelements}, it suffices to prove that there is a homotopy equivalence \[ B\left(\int_\cC F\right)\simeq B\left(\iint_\cC F\right)\] between any of the---several---models for the homotopy types of the 2-category and double category above. Indeed, this is the strategy we choose for our proof in \cref{subsec:theproof}. The next step then consists of determining which of these models lend themselves to a fruitful comparison.

\subsection{An enlightening failed attempt}\label{subsec:failed}

Before diving into the proof, we believe it is worthwhile to
share some of our failed attempts in identifying the ``correct'' models. For one, this will likely assuage some readers' curiosity, as these include some of the most natural models one may think to consider. But additionally, their shortcomings will shine a light on the desired features that the correct models should have, making our choice in \cref{subsec:theproof} seem more intuitive. Of course, any readers uninterested in this digression may skip directly to \cref{subsec:theproof}.

\begin{defn}\label{defn:2nerve}
The \textbf{bisimplicial nerve} of a 2-category $\cA$ is the bisimplicial set $N\cA$ whose set of $(m,n)$-simplices consists of pastings in $\cA$ as follows.
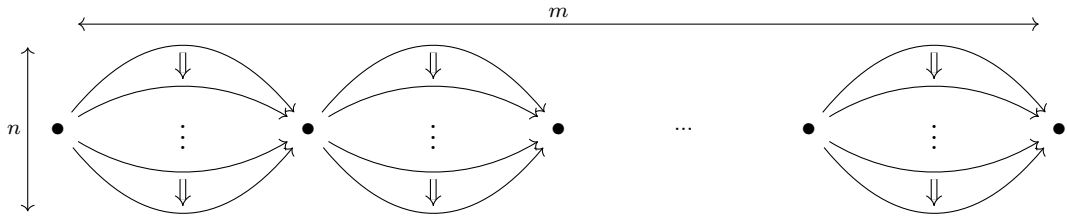
\[\begin{tikzcd}
	\bullet\ar[rrrrrrrrrrrr, shift left=15, "m", leftrightarrow] &&& \bullet &&& \bullet &&& \bullet &&& \bullet
	\arrow[""{name=0, anchor=center, inner sep=0}, curve={height=-40pt}, from=1-1, to=1-4]
	\arrow[""{name=1, anchor=center, inner sep=0}, curve={height=40pt}, from=1-1, to=1-4]
	\arrow[""{name=2, anchor=center, inner sep=0}, curve={height=-20pt}, from=1-1, to=1-4]
	\arrow[""{name=3, anchor=center, inner sep=0}, curve={height=20pt}, from=1-1, to=1-4]
	\arrow["\vdots"{description}, draw=none, from=1-1, to=1-4]
	\arrow[""{name=4, anchor=center, inner sep=0}, curve={height=-40pt}, from=1-4, to=1-7]
	\arrow[""{name=5, anchor=center, inner sep=0}, curve={height=40pt}, from=1-4, to=1-7]
	\arrow[""{name=6, anchor=center, inner sep=0}, curve={height=-20pt}, from=1-4, to=1-7]
	\arrow[""{name=7, anchor=center, inner sep=0}, curve={height=20pt}, from=1-4, to=1-7]
	\arrow["\vdots"{description}, draw=none, from=1-4, to=1-7]
	\arrow["\dots"{description}, draw=none, from=1-7, to=1-10]
	\arrow[""{name=8, anchor=center, inner sep=0}, curve={height=-40pt}, from=1-10, to=1-13]
	\arrow[""{name=9, anchor=center, inner sep=0}, curve={height=40pt}, from=1-10, to=1-13]
	\arrow[""{name=10, anchor=center, inner sep=0}, curve={height=-20pt}, from=1-10, to=1-13]
	\arrow[""{name=11, anchor=center, inner sep=0}, curve={height=20pt}, from=1-10, to=1-13]
	\arrow["\vdots"{description}, draw=none, from=1-10, to=1-13]
	\arrow[shorten <=2pt, shorten >=2pt, Rightarrow, from=0, to=2]
	\arrow[shorten <=2pt, shorten >=2pt, Rightarrow, from=3, to=1]
	\arrow[shorten <=2pt, shorten >=2pt, Rightarrow, from=7, to=5]
	\arrow[shorten <=2pt, shorten >=2pt, Rightarrow, from=4, to=6]
	\arrow[shorten <=2pt, shorten >=2pt, Rightarrow, from=8, to=10]
	\arrow[shorten <=2pt, shorten >=2pt, Rightarrow, from=11, to=9]
    \arrow[leftrightarrow, from=0, to=1, shift right=22, "n"']
\end{tikzcd}\]
\end{defn}

\begin{defn}\label{defn:dblnerve}
The \textbf{bisimplicial nerve} of a double category $\bA$ is the bisimplicial set $N\bA$ whose set of $(m,n)$-simplices consists of $m$-by-$n$ grids of squares in $\bA$.
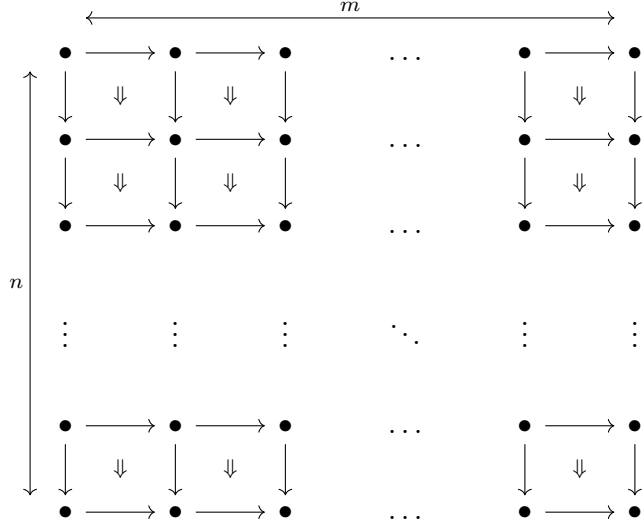
\[\begin{tikzcd}
	\bullet\ar[rrrrr, shift left=5, "m", leftrightarrow]\ar[ddddd, shift right=5, "n"', leftrightarrow]  & \bullet & \bullet & \dots & \bullet & \bullet \\
	\bullet & \bullet & \bullet & \dots & \bullet & \bullet \\
	\bullet & \bullet & \bullet & \dots & \bullet & \bullet \\
	\vdots & \vdots & \vdots & \ddots & \vdots & \vdots \\
	\bullet & \bullet & \bullet & \dots & \bullet & \bullet \\
	\bullet & \bullet & \bullet & \dots & \bullet & \bullet
	\arrow[from=1-1, to=1-2]
	\arrow[from=1-1, to=2-1]
	\arrow["\Downarrow"{description}, draw=none, from=1-1, to=2-2]
	\arrow[from=1-2, to=1-3]
	\arrow[from=1-2, to=2-2]
	\arrow["\Downarrow"{description}, draw=none, from=1-2, to=2-3]
	\arrow[from=1-3, to=2-3]
	\arrow[from=1-5, to=1-6]
	\arrow[from=1-5, to=2-5]
	\arrow["\Downarrow"{description}, draw=none, from=1-5, to=2-6]
	\arrow[from=1-6, to=2-6]
	\arrow[from=2-1, to=2-2]
	\arrow[from=2-1, to=3-1]
	\arrow["\Downarrow"{description}, draw=none, from=2-1, to=3-2]
	\arrow[from=2-2, to=2-3]
	\arrow[from=2-2, to=3-2]
	\arrow["\Downarrow"{description}, draw=none, from=2-2, to=3-3]
	\arrow[from=2-3, to=3-3]
	\arrow[from=2-5, to=2-6]
	\arrow[from=2-5, to=3-5]
	\arrow["\Downarrow"{description}, draw=none, from=2-5, to=3-6]
	\arrow[from=2-6, to=3-6]
	\arrow[from=3-1, to=3-2]
	\arrow[from=3-2, to=3-3]
	\arrow[from=3-5, to=3-6]
	\arrow[from=5-1, to=5-2]
	\arrow[from=5-1, to=6-1]
	\arrow["\Downarrow"{description}, draw=none, from=5-1, to=6-2]
	\arrow[from=5-2, to=5-3]
	\arrow[from=5-2, to=6-2]
	\arrow["\Downarrow"{description}, draw=none, from=5-2, to=6-3]
	\arrow[from=5-3, to=6-3]
	\arrow[from=5-5, to=5-6]
	\arrow[from=5-5, to=6-5]
	\arrow["\Downarrow"{description}, draw=none, from=5-5, to=6-6]
	\arrow[from=5-6, to=6-6]
	\arrow[from=6-1, to=6-2]
	\arrow[from=6-2, to=6-3]
	\arrow[from=6-5, to=6-6]
\end{tikzcd}\]
\end{defn}

Recall that, given a 2-category $\cA$, we can use the bisimplicial nerve $N\cA$ to obtain its associated topological space. In fact, there are several ways to do this: we can first realize each vertical simplicial set and then realize the resulting horizontal simplicial space, or reverse the vertical and horizontal directions in this process, or realize the diagonal simplicial set; all of these resulting spaces are homotopy equivalent. Similarly, one can use the bisimplicial nerve of a double category $N\bA$ to obtain its associated space in three ways. Hence, we could try to obtain the desired homotopy equivalence \[ B\left(\int_\cC F\right)\simeq B\left(\iint_\cC F\right)\] from some convenient bisimplicial map between $N\int_\cC F$ and $N\iint_\cC F$ in either direction. However, this does not seem feasible, as we hint at in the example below. 

First, let us introduce more formally the convenient terminology we wish to use.

\begin{notation}\label{not:horandver}
Given a bisimplicial set $X\colon\Delta^\op\times\Delta^\op\to\Set$, we think of the first simplicial direction as ``horizontal'' and the second as ``vertical'', which gives rise to horizontal face and degeneracy maps \[d_i^h\colon X_{m,n}\to X_{m-1,n}\qquad \text{ and } \qquad s_i^h\colon X_{m,n}\to X_{m+1,n},\] and similarly for their vertical counterparts \[d_i^v\colon X_{m,n}\to X_{m,n-1}\qquad \text{ and } \qquad s_i^v\colon X_{m,n}\to X_{m,n+1}.\]
\end{notation}

\begin{ex}\label{ex2}
    Let $F\colon\cC^\op\to\underline{\cat}$ be a 2-functor; we will show how the obvious attempts do not produce bisimplicial maps between $N\int_\cC F$ and $N\iint_\cC F$ in either direction.

    Unraveling the definitions, we find that both bisimplicial sets have the same set of $(0,0)$-simplices: these are the pairs $(c,x)$. Focusing on $N\int_\cC F$, we find the following:
    \begin{itemize}
        \item $(1,0)$-simplices are given by morphisms $(f,\varphi)\colon (c,x)\to (c',x')$ in $\int_\cC F$;
        \item $(0,1)$-simplices are given by objects $(c,x)$ in $\int_\cC F$;
        \item $(1,1)$-simplices are given by 2-cells in $\int_\cC F$ 
        \begin{equation}\label{eq3}
        \begin{aligned}\begin{tikzpicture}
            \node at(-0.5,0) {$(c,x)$};
            \node at(2.1,0.05) {$(c',x')$};
            \node at(0.7,0) {$\Downarrow\scriptstyle{\alpha}$};
            \draw[->] (0, 0.2) to[bend left] node[above] {$\scriptstyle{(f,\varphi)}$} (1.5,0.2);
       \draw[->] (0, -0.2) to[bend right] node[below] {$\scriptstyle{(g, (F\alpha)_{x'}\circ\varphi)}$} (1.5,-0.2);
        \end{tikzpicture}\end{aligned}
\end{equation}
        \end{itemize}
        Next, studying $N\iint_\cC F$ we find the following:
        \begin{itemize}
        \item $(1,0)$-simplices are given by horizontal morphisms $f\colon (c,Ff(x'))\to (c',x')$ in $\iint_\cC F$;
            \item $(0,1)$-simplices are given by vertical morphisms $\varphi\colon (c,x)\to (c,x')$ in $\iint_\cC F$;
        \item $(1,1)$-simplices are given by squares in $\iint_\cC F$
       \begin{equation}\label{eq4}\begin{tikzcd}
            (c,Ff(x'))\rar["f"]\dar["(F\alpha)_{y'}\circ Ff(\varphi)"']\ar[rd,phantom, "\Downarrow\scriptstyle{\alpha}" description] & (c',x')\dar["\varphi"]\\
            (c, Fg(y'))\rar["g"'] & (c', y')
            \end{tikzcd}\end{equation}
        \end{itemize}
        
If we try to construct a bisimplicial map $\Phi\colon N\int_\cC F\to N\iint_\cC F$, it seems as if the natural choice is to define the following:
\begin{itemize}
    \item $\Phi_{00}$ is the identity on $(0,0)$-simplices. 
    \item Given a $(1,0)$-simplex $(f,\varphi)\colon (c,x)\to (c',x')$, we know that $\Phi_{10}(f,\varphi)$ must be a horizontal morphism in $\iint_\cC F$; it makes sense to define it as $f\colon (c,Ff(x'))\to (c',x')$.
    \item Given a $(0,1)$-simplex $(c,x)$, we can define $\Phi_{01}(c,x)$ as the vertical morphism $\id_{(c,x)}$.
    \item Given a $(1,1)$-simplex $\alpha$ as in \cref{eq3},  the only natural square we can construct from this data is the following, as we've seen in \cref{ex1}.
    \begin{equation}\label{squareid}\begin{tikzcd}
            (c,Ff(x'))\rar["f"]\dar["(F\alpha)_{x'}"']\ar[rd,phantom, "\Downarrow\scriptstyle{\alpha}" description] & (c',x')\dar[equal]\\
            (c, Fg(x'))\rar["g"'] & (c', x')
            \end{tikzcd}\end{equation}
\end{itemize}
However, this does not give rise to a bisimplicial map as it will not respect the face maps. For instance, we have that $d_1^v \Phi_{11}(\alpha)=(F\alpha)_{x'}$, while $\Phi_{01}(d_1^v\alpha)=\id_{(c,x)}$. Similar obstructions also show that $\Phi$ does not give rise to a simplicial map between the corresponding diagonal simplicial sets: for instance, $d_0\Phi_{11}(\alpha)=(c,Ff(x'))$, but $\Phi_{00}(d_0\alpha)=(c,x)$.

What if we try to construct a map in the opposite direction, $\Theta\colon N\iint_\cC F\to \int_\cC F$? Studying the lower simplices again, the natural choice seems to be the following:
\begin{itemize}
    \item $\Theta_{00}$ is the identity on $(0,0)$-simplices.
    \item Given a $(1,0)$-simplex $f\colon (c,Ff(x'))\to (c',x')$, we can let $\Theta_{10}(f)=(f,\id)$.
    \item Given a $(0,1)$-simplex $\varphi\colon (c,x)\to (c,x')$, we have that $\Theta_{01}(\varphi)$ must be an object in $\int_\cC F$. The only natural choices are either $\Theta_{01}(\varphi)=(c,x)$ or $\Theta_{01}(\varphi)=(c,x')$; however, neither of these will give rise to a bisimplicial map. For example, using the first choice, we have that $d_0^v\Theta_{01}(\varphi)=(c,x)$, but $\Theta_{00}(d_0^v \varphi)=(c,x')$.
    \end{itemize}

Since the problem when defining $\Theta$ first came up at the level of $(0,1)$-simplices and vertical face maps, we could try to ignore these simplices and directly define $\Theta$ as a simplicial map between the diagonal simplicial sets. For this, given a 1-simplex $\alpha$ in $\diag N\iint_\cC F$ as in \cref{eq4}, we should map it to the 2-cell
\begin{center}
    \begin{tikzpicture}
            \node at(-1,0) {$(c,Ff(x'))$};
            \node at(2.1,0.05) {$(c',y')$};
            \node at(0.7,0) {$\Downarrow\scriptstyle{\alpha}$};
            \draw[->] (0, 0.2) to[bend left] node[above] {$\scriptstyle{(f,Ff(\varphi))}$} (1.5,0.2);
       \draw[->] (0, -0.2) to[bend right] node[below] {$\scriptstyle{(g, (F\alpha)_{x'}\circ Ff(\varphi))}$} (1.5,-0.2);
        \end{tikzpicture}
\end{center}
However, once again, this does not give rise to a simplicial map. To show this, consider a 2-simplex in $\diag N\iint_\cC F$ as below.
\begin{equation}\label{22cell}\begin{tikzcd}
    (c,Fgf(x''))\rar["f"]\dar["\delta"']\ar[phantom, dr, "\Downarrow\scriptstyle{\alpha}"] & (c', Fg(x''))\rar["g"]\dar["\psi"]\ar[phantom, dr, "\Downarrow\scriptstyle{\beta}"] & (c'', x'')\dar["\varphi"]\\
     (c,Fg'f'(y''))\rar["f'"]\dar[equal]\ar[phantom, dr, "\Downarrow\scriptstyle{\id}"] & (c', Fg'(y''))\rar["g'"]\dar[equal] \ar[phantom, dr, "\Downarrow\scriptstyle{\id}"] & (c'', y'')\dar[equal]\\
      (c,Fg'f'(y''))\rar["f'"'] & (c', Fg'(y''))\rar["g'"'] & (c'', y'')
\end{tikzcd}\end{equation}
In order for $\Theta$ to respect face maps, this must map to some 2-simplex in $\diag N\int_\cC F$ of the form
\[\begin{tikzcd}[column sep=small]
	{(c, Fgf(x''))} &&& {(c', Fg'(y''))} &&& {(c'', y'')}
	\arrow[""{name=0, anchor=center, inner sep=0}, "{(f', \delta)}"{description}, from=1-1, to=1-4]
	\arrow[""{name=0p, anchor=center, inner sep=0}, phantom, from=1-1, to=1-4, start anchor=center, end anchor=center]
	\arrow[""{name=1, anchor=center, inner sep=0}, "{(f, Ff(\psi))}", shift left=3, curve={height=-24pt}, from=1-1, to=1-4]
	\arrow[shift right=2, curve={height=24pt}, from=1-1, to=1-4]
	\arrow[""{name=2, anchor=center, inner sep=0}, "{(g', \id)}"{description}, from=1-4, to=1-7]
	\arrow[""{name=2p, anchor=center, inner sep=0}, phantom, from=1-4, to=1-7, start anchor=center, end anchor=center]
	\arrow[shift left=3, curve={height=-24pt}, from=1-4, to=1-7]
	\arrow[""{name=3, anchor=center, inner sep=0}, "{(g', \id)}"', shift right=2, curve={height=24pt}, from=1-4, to=1-7]
	\arrow["\alpha", shorten <=10pt, shorten >=10pt, Rightarrow, from=1, to=0p]
	\arrow["\id", shorten <=10pt, shorten >=10pt, Rightarrow, from=2p, to=3]
\end{tikzcd}\] whose total pasting is equal to the one that would result from the composite square $\alpha\vert\beta$. For that to be the case, the missing morphism on the top right must be of the form $(g, \epsilon)$ for some morphism $\epsilon\colon Fg'(y'')\to Fg(y'')$ in $c'$. However, none of the data in \cref{22cell} guarantees the existence of such a morphism.
\end{ex}

Let us draw some intuitive conclusions based on what we learned in \cref{ex1,ex2}.
\begin{enumerate}
    \item As we saw in \cref{ex1}, a morphism $(f, \varphi)\colon (c,x)\to(c',x')$ in $\int_\cC F$ can be written as the composite $(f,\varphi)=(\id,\varphi)\circ(f, \id)$ and thus encodes the same data as a ``corner''
    \[\begin{tikzcd}
        (c,x)\dar["\varphi"'] \\
        (c, Ff(x'))\rar["f"'] & (c',x')
    \end{tikzcd}\] consisting of a vertical and a horizontal morphism in $\iint_\cC F$. This discrepancy caused problems in \cref{ex2} when trying to use (the diagonal of) the bisimplicial nerves, and so it seems that we would like to have a simplicial set associated to $\int_\cC F$ whose 1-simplices are the morphisms in $\int_\cC F$, and a simplicial set associated to $\iint_\cC F$ whose 1-simplices are the corners as above.
    \item While it is true that every square in $\iint_\cC F$ gives rise to a 2-cell in $\int_\cC F$ which encodes the same data, this assignment will not respect composites in an appropriate sense, which caused problems in \cref{ex2}. Then, we expect the simplices encoding 2-cells and squares to contain some additional information that allows us to circumvent this issue.
    \item To start thinking about what that additional information might be, recall that while every 2-cell in $\int_\cC F$ gives rise to a square in $\iint_\cC F$, this process loses information due to the constraint in the objects participating in a square. Namely, if we start from a 2-cell as in \cref{eq3}, we obtain a square as in \cref{squareid}, where we have lost the data of $\varphi$ entirely. We can resolve this by adding $\varphi$ to the data that we record; that is, by considering a correspondence 
     \[\begin{tikzcd}
            (c,x)\rar[bend left, "{(f,\varphi)}"]\rar[bend right, "{(g, (F\alpha)_{x'}\circ\varphi)}"']\rar[phantom, "\Downarrow\scriptstyle{\alpha}" description] & (c',x')
        \end{tikzcd} \qquad \rightsquigarrow \qquad \begin{tikzcd}
       (c,x)\dar["\varphi"']\\
        (c,Ff(x'))\rar["f"]\dar["(F\alpha)_{x'}"']\ar[dr,phantom, "\Downarrow\scriptstyle{\alpha}" description] & (c',x')\dar[equal]\\
        (c,Fg(x'))\rar["g"'] & (c,x')
    \end{tikzcd}\]
\end{enumerate}

\subsection{Proving the main result}\label{subsec:theproof}

Guided by the intuition gained in \cref{subsec:failed}, the models for the homotopy types of $\int_\cC F$ and $\iint_\cC F$ that will allow us to prove our theorem are given by the \emph{bar construction}.  We recall this  construction below.

\begin{defn}\label{defn:bar}
The bar construction of a bisimplicial set $X$ is the simplicial set $\overline{W}X$ whose set of $k$-simplices is given by
\[(\overline{W}X)_k=\{(t_0,\dots, t_k)\in\prod_{j=0}^k X_{j,k-j} \ \vert \ d_0^v t_j=d_{j+1}^h t_{j+1} \text{ for all } 0\leq j<k\}.\] For each $0\leq i\leq k$, the face and degeneracy maps are respectively defined as
\begin{align*}
d_i(t_0,\dots,t_k) & =(d_i^v t_0, d_{i-1}^v t_1,\dots, d_1^v t_{i-1},d_i^h t_{i+1},\dots, d_i^h t_k),\\
s_i(t_0,\dots,t_k)& =(s_i^v t_0, s_{i-1}^v t_1,\dots,s_0^v t_i, s_i^h t_i,\dots, s_i^h t_k).
\end{align*}
\end{defn}

The bar construction gives us yet another possible model for the homotopy type of a bisimplicial set, as guaranteed by the next result.

\begin{theorem}\label{thm:bar}\cite[Theorem 1.1]{bar}
For any bisimplicial set $X$ there exists a natural weak homotopy equivalence $X\to\overline{W}X$.
\end{theorem}

Let us depict the lower simplices in the bar construction of the bisimplicial sets obtained by taking the diagonal of the bisimplicial nerves of a 2-category and a double category from \cref{defn:2nerve,defn:dblnerve}.

\begin{ex}\label{bar2cat}
    Let $\cA$ be a 2-category. Then, the lower simplices in $\overline{W} N\cA$ are pastings of the form below.
\[\begin{tikzcd}
 n=0: & \bullet\\
n=2: & \bullet\rar[bend right=40]\rar[bend left=40]\rar[phantom, "\Downarrow" description] & \bullet\rar & \bullet
\end{tikzcd}\hspace{1cm}\begin{tikzcd}
    n=1: & \bullet \rar & \bullet\\
n=3: & \bullet\rar\rar[phantom, bend right=35,"\Downarrow" description]\rar[phantom, bend left=35,"\Downarrow" description]\rar[bend right=75]\rar[bend left=75] & \bullet\rar[bend right=40]\rar[bend left=40]\rar[phantom, "\Downarrow" description] & \bullet\rar & \bullet
\end{tikzcd}\]
\end{ex}

\begin{ex}\label{bardoublecat}
    Let $\bA$ be a double category. Then, the  lower simplices in $\overline{W} N\bA$ are pastings of the form below.
    \[\begin{tikzcd}
         n=0: & \bullet\\
 & \bullet\dar\\
n=2: & \bullet\dar\rar\ar[dr,phantom,"\Downarrow" description]& \bullet\dar \\
&  \bullet\rar & \bullet\rar & \bullet
    \end{tikzcd}\hspace{1cm}
    \begin{tikzcd}
    & \bullet\dar \\
     n=1:  & \bullet\rar & \bullet\\
 & \bullet\dar\\
n=3: & \bullet\dar\rar\ar[dr,phantom,"\Downarrow" description] & \bullet\dar\\
& \bullet\dar\rar\ar[dr,phantom,"\Downarrow" description] & \bullet\dar\rar \ar[dr,phantom,"\Downarrow" description]& \bullet\dar\\
& \bullet\rar & \bullet\rar & \bullet\rar & \bullet
\end{tikzcd}\]
\end{ex}

The proof of \cref{thm:main} is quite technical and notation-heavy, even though the intuition behind it is likely straightforward to a reader that has followed up to this point. To illustrate this, we present a final example where we define the lower levels of inverse simplicial isomorphisms between $\overline{W} N\int_\cC F$ and $\overline{W} N \iint_\cC F$, before giving the complete proof of our main theorem with all the required verifications. 

\begin{ex}
    If we specify the lower simplices of the bar constructions described in \cref{bar2cat,bardoublecat} to the 2-category $\cA=\int_\cC F$ and double category $\bA=\iint_\cC F$, one can directly see how they encode precisely the same data, giving rise to the indicated bijective correspondences.
\[\begin{tikzcd}[column sep=small]
    n=0: & & (c,x) & \leftrightsquigarrow & (c,x)\\
    n=1: & & & & (c,x)\dar["\varphi"'] \\
    & (c,x)\rar["{(f,\varphi)}"] & (c',x') &  \leftrightsquigarrow  & (c, Ff(x'))\rar["f"] & (c',x')\\
    n=2: \\
    & &&&  (c,x)\dar["\varphi"']\\
     (c,x)\rar[bend right=40, "{(f',(F\alpha)_{x'}\circ\varphi)}"']\rar[bend left=40, "{(f, \varphi)}"]\rar[phantom, "\Downarrow" pos=0.35, "\scriptstyle{\alpha}"' pos=0.7] & (c',x')\rar["{(g, \psi)}"] & (c'',x'') & \leftrightsquigarrow & (c, Ff(x'))\rar["f"]\dar["Ff'(\psi)\circ(F\alpha)_{x'}"']\ar[rd, phantom, "\Downarrow\scriptstyle{\alpha}" description] & (c',x')\dar["\psi"]\\
     & & & & (c, Fgf'(x''))\rar["f'"'] & (c', Fg(x''))\rar["g"'] & (c'', x'')\\
     n=3:
     \end{tikzcd}\]
   \[\begin{tikzcd}[row sep=large, column sep=50]
     (c,x)\rar["{(f', (F\alpha)_{x'}\circ\varphi)}"]\rar[phantom, bend right=35,"\Downarrow" pos=0.4, "\scriptstyle{\beta}"' pos=0.6]\rar[phantom, bend left=35,"\Downarrow" pos=0.4, "\scriptstyle{\alpha}"' pos=0.6]\rar[bend right=75, "{(f'', (F\beta)_{x'}\circ (F\alpha)_{x'}\circ\varphi)}"']\rar[bend left=75, "{(f,\varphi)}"] & (c',x')\rar[bend right=40, "{(g',(F\gamma)_{x''}\circ \psi)}"']\rar[bend left=40, "{(g,\psi)}"]\rar[phantom, "\Downarrow" pos=0.4, "\scriptstyle{\gamma}"' pos=0.6] & (c'',x'')\rar["{(h,\delta)}"] & (c''',x''') \ \   \leftrightsquigarrow\end{tikzcd}\hspace{2cm}\]
     \[\hspace{2cm}\begin{tikzcd}[row sep=large, column sep=15]
     (c,x)\dar["\varphi"] \\
     (c, Ff(x'))\ar[rd, phantom, "\Downarrow\scriptstyle{\alpha}" description]\rar["f"]\dar["Ff'(\psi)\circ(F\alpha)_{x'}"'] & (c',x')\dar["\psi"]\\
     (c, Fgf'(x''))\ar[rd, phantom, "\Downarrow\scriptstyle{\beta}" description]\rar["f'"]\dar["Fg'f''(\delta)\circ (F\beta\vert\gamma)_{x''}"'] & (c', Fg(x''))\ar[rd, phantom, "\Downarrow\scriptstyle{\gamma}" description]\rar["g"]\dar["Fg'(\delta)\circ(F\gamma)_{x''}" description] & (c'', x'')\dar["\delta"]\\
     (c, Fhg'f''(x'''))\rar["f''"'] & (c', Fhg'(x'''))\rar["g'"'] & (c'', Fh(x'''))\rar["h"'] & (c''',x''')\\
\end{tikzcd}\]

\end{ex}

\begin{proof}[Proof of \cref{thm:main}]
We will prove that there is a weak homotopy equivalence $B\left(\int_\cC F\right)\simeq B\left(\iint_\cC F\right)$ by constructing an isomorphism of simplicial sets \[\Phi\colon \overline{W} N\int_\cC F \rightleftarrows \overline{W} N\iint_\cC F\colon \Theta\] and appealing to \cref{thm:bar}.

We start by explicitly describing the $p$-simplices of both bar constructions, and introducing our notation. Unraveling \cref{defn:2catelt,defn:2nerve,defn:bar}, we see that a $p$-simplex of $\overline W N \int_\cC F$ consists of the data
\begin{equation}\label{simpl2cat}
    \left\{\begin{array}{ll}
        (c_m, x_m) & 0 \leq m \leq p\\
        (f^n_m, \varphi^n_m)\colon (c_m, x_m) \to (c_{m-1}, x_{m-1}) & 1 \leq m \leq p, 0 \leq n < m\\
        \alpha^n_m\colon (f^{n+1}_m, \varphi^n_m) \Rightarrow (f^{n}_m, \varphi^{n-1}_m) & 2 \leq m \leq p, 0 \leq n < m-1
    \end{array}\right.
\end{equation}
where each $(c_m, x_m)$, $(f^n_m, \varphi^n_m)$, and $\alpha^n_m$ are, respectively, objects, morphisms, and 2-cells in $\int_\cC F$. The components in these data are subject to the relations given by the 2-cells $\alpha_m^n$ as specified in \cref{defn:2catelt}. Namely, we must have that \[\varphi_m^{n}=(F\alpha_m^{n})_{x_{m-1}}\circ \varphi_m^{n+1}\] for all $1\leq m \leq p$ and all $0\leq n\leq m-2$. Hence, iterating this relation we obtain that the second component of any morphism $(f^n_m, \varphi^n_m)$ in \cref{simpl2cat} for $0\leq n\leq m-2$ can be determined from the corresponding ``top'' morphism $\varphi_m^{m-1}$ together with all the 2-cells $\alpha_m^i$, as we have that \begin{equation}\label{eq:2reduced}
(f^n_m, \varphi^n_m) = \left(f^n_m, (F\alpha_m^{n}\circ\dots \circ \alpha^{m-2}_m)_{x_{m-1}}\circ \varphi^{m-1}_m\right).\end{equation}

Similarly, unraveling \cref{defn:dblcatelt,defn:dblnerve,defn:bar}, we see that a $p$-simplex of $\overline W N \iint_\cC F$ consists of the following information
\begin{equation}\label{simpldblcat}
    \left\{\begin{array}{ll}
        (c_m, x^n_m) & 0 \leq m \leq p, 0 \leq n \leq m\\
        f^n_m\colon (c_m, x^n_m) \to (c_{m-1}, x^n_{m-1}) & 1 \leq m \leq p, 0 \leq n < m\\
        \varphi^n_m\colon (c_m, x^{n+1}_m) \to (c_m, x^{n}_m) & 1 \leq m \leq p, 0 \leq n < m\\
        \alpha^n_m\colon (\varphi^{n}_m \; ^{f^{n+1}_m}_{\substack{f^{n}_m}} \; \varphi^{n}_{m-1}) & 2 \leq m \leq p, 0 \leq n < m-1
    \end{array}\right.
\end{equation} 
where each $(c_m, x^n_m)$, $f^n_m$, $\varphi^n_m$, and $\alpha^n_m$ are, respectively, objects, horizontal morphisms, vertical morphisms, and squares in $\iint_\cC F$. Once again, these components are subject to additional relations imposed by the structure in the double category $\iint_\cC F$. First, due to the definition of the horizontal morphisms in the double category, we know that the objects must satisfy \[x_m^n=Ff_m^n(x_{m-1}^n)\] for all $1\leq m\leq p$ and all $0\leq n\leq m-1$. Thus, iterating this we obtain that any object $(c_m, x_m^n)$ can be determined from the corresponding ``top'' object $(c_n, x_n^n)$ together with all the morphisms $f_i^n$, due to the equation
\begin{equation}\label{eq:objreduced}
    (c_m,x_m^n)=(c_m, F f^n_{n+1}\dots f_m^n (x_n^n)).
\end{equation}
Moreover, the conditions imposed by the squares $\alpha_m^n$ give us that \[\varphi_m^n=Ff_{m}^{n}(\varphi_{m-1}^n)\circ (F\alpha_m^{n})_{x_{m-1}^{n+1}}\] for all $2\leq m\leq p$ and $0\leq n\leq m-2$. Iterating, we can see that any such vertical morphisms $\varphi_m^n$ may be determined from the corresponding ``rightmost'' vertical morphism $\varphi_{n+1}^n$ together with all the horizontal morphisms $f_i^{n-1}$ and all squares $\alpha_i^n$, due to the equation
\begin{equation}\label{eq:vermorreduced}
    \varphi_m^n=Ff_{n+2}^{n}\dots f_m^{n} (\varphi_{n+1}^n)\circ (F\alpha_m^n\vert\dots\vert\alpha_{n+2}^n)_{x_{n+1}^{n+1}}.
\end{equation}

We now proceed to construct the bijection. For each $p$, define a function \[\Phi_p\colon \left(\overline W N \int_\cC F\right)_p \to \left(\overline W N \iint_\cC F\right)_p\] mapping a $p$-simplex as in \cref{simpl2cat} to
\begin{equation}
    \begin{cases}
       \begin{cases}
       (c_m, x_m) \\
       (c_m, Ff^n_{n+1}\dots f^n_m(x_n)) 
       \end{cases} & \begin{aligned} & 0  \leq m\leq p, n=m \\
 & 0  \leq m \leq p, 0 \leq n < m  \end{aligned}\\
        f^n_m & 1 \leq m \leq p, 0 \leq n < m\\
        \begin{cases}
        \varphi_m^{m-1} \\
        Ff^{n}_{n+2}\dots f^{n}_m(\varphi^{n}_{n+1})\circ (F\alpha^n_m\vert\dots\vert\alpha^n_{n+2})_{x_{n+1}} \end{cases} & \begin{aligned} & 1 \leq m \leq p, n=m-1\\
  & 1 \leq m \leq p, 0\leq n < m-1 \end{aligned}\\ 
        \alpha^n_m &2 \leq m \leq p, 0 \leq n < m-1
   \end{cases}
\end{equation} 
In the opposite direction, consider the function \[\Theta_p\colon \left(\overline W N \iint_\cC F\right)_p \to \left(\overline W N \int_\cC F\right)_p\] mapping a $p$-simplex as in \cref{simpldblcat} to
\begin{equation}\label{theta}
   \begin{cases}
        (c_m, x^m_m) & 0\leq m \leq p\\
        \begin{cases}
            (f_m^{m-1}, \varphi_m^{m-1})\\
    (f^n_m, (F\alpha^{n}_m \circ \dots \circ \alpha^{m-2}_m)_{x^{m-1}_{m-1}} \circ \varphi^{m-1}_m) \end{cases} & \begin{aligned}
        & 1\leq m\leq p, n=m-1\\
& 1\leq m \leq p, 0\leq n < m-1
    \end{aligned}\\
        \alpha^n_m & 2\leq m \leq p, 0\leq n < m-1
  \end{cases}
\end{equation}

We claim that $\Phi_p$ and $\Theta_p$ are inverses for every $p$. To see this, first examine the composite $\Theta_p \Phi_p$. Applying this composite to a $p$-simplex as in \cref{simpl2cat}, we obtain a $p$-simplex with the same data at the level of objects and of 2-cells as the initial one. To check that the data corresponding to morphisms also agrees, note that the first components are always of the form $f_m^n$ and thus agree with the initial data. It then suffices to show that we recover the same ``top'' morphisms $(f_m^{m-1}, \varphi_m^{m-1})$ for all $m$, as the relations in \cref{eq:2reduced} would then imply that all other morphism data must also agree. It is straightforward to verify that this is indeed the case, and thus $\Theta_p\Phi_p=\id$ for all $p$.

We now focus on the composite $\Phi_p \Theta_p$. Applying this composite to a $p$-simplex as in \cref{simpldblcat}, we get a $p$-simplex with the same data at the level of horizontal morphisms and squares. Additionally, one can easily see that the data of the ``top'' objects $(c_m, x_m^m)$ also agrees for every $m$; hence, the data of all objects $(c_m, x_m^n)$ must also agree, as these are determined from the ``top'' objects by \cref{eq:objreduced}. Similarly, it is not hard to see that the data of the ``rightmost'' vertical morphisms $\varphi_m^{m-1}$ agrees for all $m$; then, the data of all vertical morphisms must also agree, since these are determined from the data recovered so far by \cref{eq:vermorreduced}. This proves that $\Phi_p\Theta_p=\id$ for all $p$, and concludes the verification that $\Phi_p$ and $\Theta_p$ are inverses of each other.

To conclude the proof, it remains to show that the maps $\Phi$ and $\Theta$ are simplicial; that is, they commute with all face and degeneracy maps. Since we have already established that they are inverses levelwise, it suffices to show that one of them is simplicial, as this implies that the other one must be as well---because the inverses of the components of a natural isomorphism always assemble into a natural isomorphism. We will prove that $\Theta$ is a simplicial map.

Beginning with face maps, consider a $p$-simplex as in \cref{simpldblcat}. If we first apply the face maps $d_i$ described in \cref{defn:bar}, we get that the resulting $(p-1)$-simplex consists of the following data when $i=0,p$,
\[\begin{cases}
    (c_m,x_m^n) & m,n\neq p-i\\
    f_m^n & m,n\neq p-i\\
    \varphi_m^n & m,n \neq p-i\\
    \alpha_m^n & m,n\neq p-i
\end{cases}
\]
and the data below when $0<i<p$.
\[\begin{cases}
    (c_m, x_m^n) & m,n\neq p-i\\
    \begin{cases}
        f_m^n\\
        f_{p-i}^n f_{p-i+1}^n
    \end{cases} & \begin{aligned}
       & m\neq p-i, p-i+1, n\neq p-i\\
       & \text{else}
    \end{aligned}\\
    \begin{cases}
        \varphi_m^n\\
        \varphi_m^{p-i-1}\varphi_m^{p-i}
    \end{cases} & \begin{aligned}
        & m\neq p-i, n\neq p-i-1, p-i\\
        & \text{else}
    \end{aligned} \\
    \begin{cases}
        \alpha_m^n\\
        \alpha_m^{p-i-1}\circ\alpha_m^{p-i}\\
        \alpha_{p-i+1}^n\vert\alpha_{p-i}^n
    \end{cases} & \begin{aligned}
        & m\neq p-i, p-i+1, n\neq p-i-1, p-i\\
        & p-i+1 < m\leq p\\
        & \text{else}
    \end{aligned}
\end{cases}\]
For ease of notation, we have omitted from the notation above the usual ranges of $m$ and $n$, only highlighting the distinguishing features of each case. Applying the map $\Theta_{p-1}$ to these simplices, we obtain, for the first case,
\begin{equation}\label{d0}\begin{cases}
    (c_m, x_m^m) & m\neq p-i\\
    \begin{cases}
        (f_m^{m-1}, \varphi_m^{m-1}) \\
        (f_m^n, -)
    \end{cases} & \begin{aligned}
        & m,n \neq p-i, n=m-1\\
        & m,n\neq p-i, n< m-1
    \end{aligned} \\
    \alpha_m^n & m,n\neq p-i
\end{cases}\end{equation}
where we have used a dash in place of the redundant information that is uniquely determined by the rest of the data. For the second case, applying $\Theta_{p-1}$ yields 
\begin{equation}\label{di}\begin{cases}
    (c_m, x_m^m) & m,n\neq p-i\\
    \begin{cases}
        (f_m^{m-1}, \varphi_m^{m-1})\\
        (f_{p-i}^{p-i-1}f_{p-i+1}^{p-i-1}, \varphi_{p-i+1}^{p-i-1}\varphi_{p-i+1}^{p-i})\\
        (f_m^{m-1}, -)\\
        (f_{p-i}^{p-i-1}f_{p-i+1}^{p-i-1}, -)
    \end{cases} & \begin{aligned}
        & m\neq p-i, p-i+1, n=m-1\\
        & m=p-i+1, n=m-1\\
        & m\neq p-i, p-i+1, n\neq p-i, n<m-1\\
        & m=p-i+1, n<m-1\\
    \end{aligned}\\
     \begin{cases}
        \alpha_m^n\\
        \alpha_m^{p-i-1}\circ\alpha_m^{p-i}\\
        \alpha_{p-i+1}^n\vert\alpha_{p-i}^n
    \end{cases} & \begin{aligned}
        & m\neq p-i, p-i+1, n\neq p-i-1, p-i\\
        & p-i+1 < m\leq p\\
        & \text{else}
    \end{aligned}
\end{cases}\end{equation}

On the other hand, if we start with a $p$-simplex as in \cref{simpldblcat} and first apply $\Theta_p$, we obtain a $p$-simplex as in \cref{theta}. It is not hard to verify that applying the face maps $d_0$ and $d_p$ to this $p$-simplex yields precisely the same $(p-1)$-simplex described in \cref{d0}. If we instead apply an inner face map $d_i$ to the $p$-simplex in \cref{theta}, the resulting $(p-1)$-simplex consists of the following data:
\[\begin{cases}
    (c_m, x_m) & m\neq p-i\\
    \begin{cases}
        (f_m^{m-1}, \varphi_m^{m-1})\\
        (f_{p-i}^{p-i-1}, \varphi_{p-i}^{p-i-1})\circ(f_{p-i+1}^{p-i-1},(F\alpha_{p-i+1}^{p-i-1})_{x_{p-i}^{p-i}}\circ\varphi_{p-i+1}^{p-i})\\
        (f_m^{m-1}, -)\\
        (f_{p-i}^{p-i-1}, -)\circ (f_{p-i+1}^{p-i-1}, -)
    \end{cases} & \begin{aligned}
        & m\neq p-i, p-i+1, n=m-1\\
        & m=p-i+1, n=m-1\\
        & m\neq p-i, p-i+1, n\neq p-i, n<m-1\\
        & m=p-i+1, n<m-1\\
    \end{aligned}\\
    \begin{cases}
        \alpha_m^n\\
        \alpha_m^{p-i-1}\circ\alpha_m^{p-i}\\
        \alpha_{p-i+1}^n\vert\alpha_{p-i}^n
    \end{cases} & \begin{aligned}
        & m\neq p-i, p-i+1, n\neq p-i-1, p-i\\
        & p-i+1 < m\leq p\\
        & \text{else}
    \end{aligned}
\end{cases}\]

We can immediately see that this $(p-1)$-simplex and the one from \cref{di} contain the same data at the level of objects and 2-cells, and that the maps $f_j^i$ in the first coordinate of all morphisms also agree. Then, by \cref{eq:2reduced}, it suffices to show that the ``top'' second coordinates in the morphisms agree as well; that is, we must compare these morphisms in the cases when $n=m-1$. According to our descriptions, these are identical when $m\neq p-i, p-i+1$. For the case $m=p-i+1$, the following equations show that they agree as well:
\begin{align*}
    (f_{p-i}^{p-i-1}, \varphi_{p-i}^{p-i-1})\circ(f_{p-i+1}^{p-i-1}, &(F\alpha_{p-i+1}^{p-i-1})_{x_{p-i}^{p-i}}\circ\varphi_{p-i+1}^{p-i}) \\ 
    & = (f_{p-i}^{p-i-1}f_{p-i+1}^{p-i-1}, Ff_{p-i+1}^{p-i-1}(\varphi_{p-i}^{p-i-1})\circ (F\alpha_{p-i+1}^{p-i-1})_{x_{p-i}^{p-i}}\circ \varphi_{p-i+1}^{p-i})\\
    & = (f_{p-i}^{p-i-1}f_{p-i+1}^{p-i-1}, \varphi_{p-i+1}^{p-i-1}\circ\varphi_{p-i+1}^{p-i}).
\end{align*}
Here, the first equality is using the definition of the composition of morphisms in the 2-category of elements, and the second equality is imposed by the relation given by the square $\alpha_{p-i+1}^{p-i-1}$. This concludes the proof that $\Theta$ respects face maps.

Now turning to degeneracy maps, again consider the $p$-simplex from  \cref{simpldblcat}. If we first apply the degeneracy maps $s_i$, the resulting $(p+1)$-simplex consists of the following data:
\[\begin{cases}
    (c_m, x_m^n)  \text{ with all } (c_{p-i}, x_{p-i}^n) \text{ and } (c_m, x_m^{p-i}) \text{ repeated}\\
    \text{all } \id_{c_{p-i}} \text{ and } f_m^n, \text{ with all } f_m^{p-i} \text{ repeated} \\
    \text{all } e_{x_m^{p-i}} \text{ and } \varphi_m^n, \text{ with all } \varphi_{p-i}^n \text{ repeated} \\
    \text{all } e_{f_m^{p-i}}, \id_{\varphi_{p-i}^n}, \text{ and } \alpha_m^n
\end{cases}\]
Applying $\Theta_{p+1}$ results in the simplex given by the following data:
\begin{equation}\label{deg1}\begin{cases}
    (c_m, x_m^m)  \text{ with all } (c_{p-i}, x_{p-i}^{p-i}) \text{ repeated}\\
    \begin{cases}
        \text{all } (f_m^{m-1}, \varphi_m^{m-1}) \text{ and } (\id_{c_{p-i}}, e_{x_{p-i}^{p-i}}) \text{ along the top}\\ 
        \text{all } (f_m^n,-) \text{ and } (\id_{c_{p-i}},-) \text{ below the top, with } (f_m^{p-i},-) \text{ repeated}
    \end{cases}\\
    \text{all } e_{f_m^{p-i}}, \id_{\varphi_{p-i}^n}, \text{ and } \alpha_m^n
\end{cases}\end{equation}

On the other hand, we can first apply $\Theta_{p}$ to a $p$-simplex as in \cref{simpldblcat} to obtain a $p$-simplex as in \cref{theta}. Now, applying the degeneracy map $s_i$ yields a $(p+1)$-simplex whose data is given as follows:
\begin{equation}\label{deg2}\begin{cases}
    (c_m, x_m^m)  \text{ with all } (c_{p-i}, x_{p-i}^{p-i}) \text{ repeated}\\
    \begin{cases}
        \text{all } (f_m^{m-1}, \varphi_m^{m-1}) \text{ and } \id_{(c_{p-i}, x_{p-i}^{p-i})} \text{ along the top}\\ 
        \text{all } (f_m^n,-) \text{ and } \id_{(c_{p-i}, x_{p-i}^{p-i})} \text{ below the top, with } (f_m^{p-i},-) \text{ repeated}
    \end{cases}\\
    \text{all } e_{f_m^{p-i}}, \id_{\varphi_{p-i}^n}, \text{ and } \alpha_m^n
\end{cases}\end{equation}
By definition of the identity morphisms in the 2-category of elements, we have that \[(\id_{c_{p-i}}, e_{x_{p-i}^{p-i}})=\id_{(c_{p-i}, x_{p-i}^{p-i})}\] from which we see that the simplices in \cref{deg1,deg2} agree. This implies that $\Theta$ respects degeneracies, which shows it is a simplicial isomorphism, concluding our proof. 
\end{proof}

\bibliographystyle{alpha}
\bibliography{references}

\end{document}